\documentclass[12 pt,twoside]{article}
\usepackage{graphicx}
\newtheorem{theorem}{\bf Theorem}[section]

\newtheorem{lemma}[theorem]{\bf Lemma}
\newtheorem{corollary}[theorem]{\bf Corollary}

\newtheorem{conjecture}{\bf Conjecture}

\input{amssym}

\begin{document}

\title{{\normalsize \bf On group choosability of total graphs}}
\bigskip\author{\small H.J. Lai$^{\textrm{a}}$, G.R. Omidi$^{\textrm{b},\textrm{c},{1}}$, G. Raeisi$^{\textrm{b}}$, \\
\footnotesize  $^{\textrm{a}}$ Department of Mathematics, West
Virginia
University,\\ \footnotesize  Morgantown, WV 26505, USA\\
\footnotesize  $^{\textrm{b}}$Department of Mathematical Sciences,
Isfahan
University of Technology,\\ \footnotesize Isfahan, 84156-83111, Iran\\
\footnotesize  $^{\textrm{c}}$School of Mathematics, Institute for
Research
in Fundamental Sciences (IPM),\\
\footnotesize  P.O.Box: 19395-5746, Tehran, Iran\\
\footnotesize {hjlai@math.wvu.edu, romidi@cc.iut.ac.ir,
g.raeisi@math.iut.ac.ir}}
\date {}
\footnotesize\maketitle\footnotetext[1] {This research was in part
supported by a grant from IPM (No.89050037)} \vspace*{-0.5cm}

\footnotesize
\begin{abstract}\rm{}

\footnotesize In this paper, we study the group and list group
colorings of total graphs and we give two group versions of the
total and list total colorings conjectures.  We establish the
group version of the total coloring conjecture for the following classes of
graphs: graphs with small maximum degree, two-degenerate graphs,
planner graphs with maximum degree at least 11, planner graphs
without certain small cycles, outerplanar and near-outerplanar
graphs. In addition, the group version of the list total
coloring conjecture is established for forests, outerplanar graphs and graphs with maximum degree at most two.


\end{abstract}
\small
\section{\large{Introduction}}
\medskip

\hspace{.5 cm}Throughout this paper, we consider simple graphs,
the graphs without loops and multiple edges and we follow
\cite{diestel} for terminology and notation not defined here. For
a graph $G$, we denote its vertex set, edge set, maximum degree
and minimum degree by $V(G)$, $E(G)$, $\Delta(G)$ and $\delta(G)$ (or simply $V$, $E$, $\Delta$, $\delta$),
respectively. If $v\in V(G)$, we use $\deg_G{(v)}$ (or simply
$\deg{(v)}$) and $N_G(v)$ to denote the degree and neighborhoods
of $v$ in $G$, respectively.

\bigskip A {\it proper coloring} of $G$ is a coloring of the
vertices of $G$ so that no two adjacent vertices are assigned the
same color. The minimum number of colors in any proper coloring of
$G$, $\chi(G)$, is called the {\it chromatic number} of $G$. A {\it
$k$-list assignment} for a graph $G$ is a function $L$ which assigns
to each vertex $v$ of $G$ a list of colors $L(v)$ such that
$|L(v)|=k$. An {\it $L$-coloring} of $G$ is a proper coloring $c$
such that $c(v)\in L(v)$ for each vertex $v$. If for every $k$-list
assignment $L$, a proper $L$-coloring of $G$ exists, then $G$ is
said to be {\it $k$-choosable} and the {\it choice number},
$\chi_l(G)$, of $G$ is the smallest $k$ such that $G$ is
$k$-choosable. The concepts of {\it chromatic index}, $\chi^\prime(G)$, and
{\it choice index}, $\chi_l^\prime(G)$, can be defined similarly in
terms of coloring the edges of $G$.

\bigskip Recall that the {\it total graph} of a graph $G$, denoted by
$T(G)$, is a graph where its vertices are the edges and vertices of
$G$ and adjacency in $T(G)$ is defined as adjacency or incidence for
the corresponding elements of $G$. The {\it total chromatic number}
of $G$, $\chi^{\prime\prime}(G)$, is the chromatic number of $T(G)$.
Clearly $\chi^{\prime\prime}(G)=\chi(T(G))\geq \Delta(G)+1$. Behzad
\cite{behzad} and Vizing \cite{Vizing} posed independently the
following famous conjecture, which is known as the {\it total coloring
conjecture}.

\begin{conjecture}\rm{}
For any graph $G$, $\chi^{\prime\prime}(G)\leq \Delta(G)+2$.
\end{conjecture}

The {\it total choice number} of $G$, $\chi_l^{\prime\prime}(G)$, is
the choice number of $T(G)$. It follows directly from the definition
that $\chi_l^{\prime\prime}(G)\geq\chi^{\prime\prime}(G)$. The
notation of total choosability was first introduced by Borodin et
al. \cite{list total conj}. They proposed the following conjecture,
known as the {\it list total coloring conjecture}.

\begin{conjecture}\rm{}
For any graph $G$, $\chi_l^{\prime\prime}(G)=\chi^{\prime\prime}(G)$.
\end{conjecture}

The concept of group coloring of graphs was first introduced by
Jaeger et al. \cite{group}. Assume that $A$ is a group and $F(G,A)$
denotes the set of all functions $f:E(G)\longrightarrow A$. Consider
an arbitrary orientation of $G$. Graph $G$ is called {\it
$A$-colorable} if for every $f\in F(G,A)$, there is a vertex
coloring $c:V(G)\longrightarrow A$ such that $c(x)-c(y)\neq f(xy)$
for each directed edge from $x$ to $y$. The {\it group chromatic
number} of $G$, $\chi_g(G)$, is the minimum $k$ such that $G$ is
$A$-colorable for any group $A$ of order at least $k$. In
\cite{group choosablity}, the concept of  group choosablity is
introduced as an extension of list coloring and group coloring. Let
$A$ be a group of order at least $k$ and $L:V(G)\longrightarrow 2^A$
be a list assignment of $G$. For $f\in F(G,A)$, an {\it
$(A,L,f)$-coloring} of $G$ is an $L$-coloring $c:V(G)\longrightarrow
A$ such that $c(x)-c(y)\neq f(xy)$ for each directed edge from $x$
to $y$. If for each $f\in F(G,A)$ there exists an $(A,L,f)$-coloring
for $G$, then we say that $G$ is {\it $(A,L)$-colorable}. If for any
group $A$ of order at least $k$ and any $k$-list assignment
$L:V(G)\longrightarrow 2^A$, $G$ is $(A,L)$-colorable, then we say
that $G$ is {\it $k$-group choosable}. The {\it group choice number}
of $G$, $\chi_{gl}(G)$, is the smallest $k$ such that $G$ is
$k$-group choosable. Clearly group choosblity of a graph is
independent of the orientation on $G$. The concept of group
choosability is also studied in \cite{On group choosability of
graphs}. The authors used the concept of {\it $D$-group
choosability} to establish the group version of the Brook's Theorem.
A graph $G$ is called {\it $D$-group choosable} if it is
$(A,L)$-colorable for each group $A$ with $|A|\geq\Delta(G)$ and
every list assignment $L:V(G)\longrightarrow 2^A$ with
$|L(v)|=\deg(v)$. They proved the following theorem, which is a
characterization of $D$-group choosable graphs.

\begin{theorem}\rm{}\label{D-choosable}(\cite{On group choosability of graphs})
The connected graph $G$ is $D$-group choosable if and only if $G$ has a block which is neither a complete graph nor a cycle.
\end{theorem}

 The following result is the group version of the Brook's Theorem.
\begin{theorem}\rm{}(\cite{On group choosability of graphs})\label{Brook's theorem}
For every connected simple graph $G$, $\chi_{gl}(G)\leq\Delta(G)+1,$ with equality if and only if $G$ is either a cycle or a complete graph.
\end{theorem}

We extend the concepts of total and list total colorings to total
group and list total group colorings of graphs. We let
$\chi_{g}^{\prime\prime}(G)=\chi_{g}(T(G))$ (resp.
$\chi_{gl}^{\prime\prime}(G)=\chi_{gl}(T(G))$) and we called it
the {\it total group chromatic number} (resp. {\it total group
choice number}) of $G$. Clearly the following inequality holds for
the mentioned chromatic numbers of $G$.
$$\chi_{gl}^{\prime\prime}(G)\geq \max\{\chi_{g}^{\prime\prime}(G),\chi_{l}^{\prime\prime}(G)\}\geq \chi^{\prime\prime}(G).$$

Now we extend the total coloring and list total coloring
conjectures as follows:
\medskip

\begin{conjecture}\rm{}\label{total conjecture 1}
For every graph $G$, $\chi_{gl}^{\prime\prime}(G)\leq
\Delta(G)+2$.
\end{conjecture}
\begin{conjecture}\rm{}\label{total conjecture 2}
For every graph $G$, $\chi_{gl}^{\prime\prime}(G)=\chi_{g}^{\prime\prime}(G)$.
\end{conjecture}

The following conjecture express the weaker version of Conjecture \ref{total conjecture 1}.
\medskip
\begin{conjecture}\rm{}\label{total conjecture 3}
For every graph $G$, $\chi_{g}^{\prime\prime}(G)\leq \chi^{\prime\prime}(G)+1$.
\end{conjecture}

In this paper, we are interested in Conjectures \ref{total
conjecture 1} and \ref{total conjecture 2} and we will establish
conjecture \ref{total conjecture 1} for certain classes of graphs
such as planar graphs with maximum degree at least 11,
two-degenerate graphs, planar graphs without certain cycles,
outerplanar and near-outerplanar graphs. Also we show that
Conjecture \ref{total conjecture 2} holds for graphs with maximum
degree at most two, forests and outerplanar graphs. Subsequently, it
will be shown that Conjecture \ref{total conjecture 3} holds for the
above mentioned classes of graphs.

\section{\large{Some upper bounds}}
\medskip

\hspace{.5 cm}In this section, we give some upper bounds for
$\chi_{gl}^{\prime\prime}(G)$ of a graph $G$ and use these bounds
to verify Conjectures \ref{total conjecture 1} and \ref{total
conjecture 2} for some classes of graphs. The {\it group choice
index} of a graph $G$, $\chi_{gl}^\prime(G)$, is defined as the group
choice number of its line graph i.e.
$\chi_{gl}^\prime(G)=\chi_{gl}(L(G))$. Clearly
$\chi_{gl}^\prime(G)\geq \chi^\prime(G)\geq \Delta(G)$. This concept
is studied in \cite{edge-group choosablity} where the authors
conjectured that every graph with maximum degree $\Delta$ is
$(\Delta+1)$-edge group choosable. Moreover, they gave infinite
families of graphs $G$ with $\chi_{gl}^\prime(G)=\Delta(G)$.
The following lemma shows that Conjecture \ref{total conjecture 1} holds for
these graphs.

\begin{lemma}\rm{}\label{upper bound1}
For every graph $G$ we have $\chi_{gl}^{\prime\prime}(G)\leq \chi_{gl}^{\prime}(G)+2$.
\end{lemma}
\textbf{Proof. }Let $G$ be a graph and $A$ be a group of order at
least $\chi_{gl}^{\prime}(G)+2$. Also let $L:V\cup E\longrightarrow
2^A$ be any ($\chi_{gl}^{\prime}(G)+2$)-list assignment of $V\cup E$
and $f\in F(T(G),A)$ be arbitrary. First we color the vertices of
$G$ from their lists. Since $\chi_{gl}^{\prime}(G)\geq \Delta(G)$
and $\chi_{gl}(G)\leq \Delta(G)+1$ by Theorem \ref{Brook's theorem},
such a coloring exists. Let $c:V(G)\longrightarrow A$ be such a
coloring so that $c(v)\in L(v)$. For each edge $e=uv$ of $G$,
without loss of generality, let the edge $eu$ be directed from $e$
to $u$ and also the edge $ev$ be directed from $e$ to $v$ in $T(G)$.
For each edge $e=uv$ of $G$, remove $f(eu)+c(u)$ and $f(ev)+c(v)$
from $L(e)$. Since for each edge $e$ of $G$, $|L(e)|\geq
\chi_{gl}^{\prime}(G)+2$, each edge of $G$ retains at least
$\chi_{gl}^{\prime}(G)$ admissible colors in its list and so, by the
definition of $\chi_{gl}^{\prime}(G)$, it is possible to color the
edges of $G$ from their lists. So we can color the vertices of
$T(G)$ from their lists and this yields an $(A,L,f)$-coloring of
$T(G)$, which shows that $\chi_{gl}^{\prime\prime}(G)\leq
\chi_{gl}^{\prime}(G)+2$.

$\hfill \blacksquare$

\medskip
The {\it coloring number} of $G$, $\mbox{col}(G)$, is the smallest
integer for which there exists an ordering of the vertices of $G$
such that each vertex $v$ has at most $\mbox{col}(G)-1$ neighbors
among vertices of smaller indices. A graph $G$ is called {\it
$d$-degenerate} if $d\geq\mbox{col}(G)-1$. For a non-regular graph
$G$ with maximum degree $\Delta$, it is easy to see that
$\mbox{col}(T(G))\leq \Delta+\mbox{col}(G)-1.$ If $G$ is a
$\Delta$-regular graph, then $T(G)$ is $2\Delta$-regular which is
not complete or a cycle since $|V(G)|\geq 4$. So by Theorem
\ref{Brook's theorem}, we have $\chi_{gl}^{\prime\prime}(G)\leq
2\Delta=\Delta+\mbox{col}(G)-1.$ Hence we have the following result.

\medskip
\begin{theorem}\rm{}\label{bound in term of col(G)}
Let $G$ be a graph with maximum degree $\Delta$. Then $$\chi_{gl}^{\prime\prime}(G)\leq \Delta+\mbox{col}(G)-1.$$
\end{theorem}

\medskip
Our first application of Theorem \ref{bound in term of col(G)}
is the following, which states that Conjectures \ref{total
conjecture 1} and \ref{total conjecture 2} hold for all forests.

\medskip
\begin{corollary}\rm{}\label{total of tree}
Let $G$ be a forest with maximum degree $\Delta\geq 2$. Then
$$\chi_{gl}^{\prime\prime}(G)=\chi_{g}^{\prime\prime}(G)=\Delta+1.$$
\end{corollary}

\bigskip
A graph $H$ is a {\it minor} of a graph $K$ if $H$ can be obtained
from a subgraph of $K$ by contracting some edges. A graph $G$ is
called {\it $K_4$-minor free} if it has no subgraph isomorphic to a
minor of $K_4$. It is well-known that \cite{series-parallel
networks} every  $K_4$-minor free graph has a vertex of degree at
most two.  A planar graph is called {\it{outerplanar}} if it has a
drawing in which each vertex lies on the boundary of the outer face.
It is well-known that a graph is outerplanar if and only if it
contains neither $K_4$ nor $K_{2, 3}$ as minors (see for example
\cite{diestel}).

\bigskip
It is well known that  \cite{2-choosable} a connected graph $G$ is
2-choosable if and only if $G$ obtained by successively removing
vertices of degree 1 until what remains, is isomorphic to either
$K_1$, $C_{2m+2}$ or $\Theta_{2,2,2m}$ for some $m$, where
$\Theta_{2,2,2m}$ is the graph consisting of two distinguished
vertices $v_0,v_{2m}$ connected by three paths $P_1$, $P_2$ and
$P_3$ of lengths 2, 2 and $2m$, respectively. So the class of
two-degenerate graphs properly contain, 2-choosable graphs,
outerplanar graphs, non-regular subcubic graphs, planar graphs of
girth at least six and {\it unicycle graphs}, graphs with exactly
one cycle, as subclasses. Using Theorem \ref{bound in term of
col(G)}, we have the following corollary.

\medskip
\begin{corollary}\rm{}\label{col(G)<3}
Conjecture \ref{total conjecture 1} holds for every two-degenerate graph. In particular, Conjecture \ref{total
conjecture 1} holds for planar graphs of girth at least six, $K_4$-minor free, outerplanar and 2-choosable graphs.
\end{corollary}
\section{\large{Graphs with bounded degrees}}
\bigskip

\hspace{.5 cm}In this section, we prove that Conjectures \ref{total conjecture 1}
and \ref{total conjecture 2} hold for graphs with maximum degree at
most two and wheel graphs with maximum degree at least five. Also
we prove that Conjecture \ref{total conjecture 1} holds for
any planar graph with maximum degree at least 11.

\begin{theorem}\rm{}\label{total of C_n}
Let $P_n$ and $C_n$ be path and cycle on $n\geq 2$
vertices, respectively. Then
$\chi_g^{\prime\prime}(P_n)=\chi_{gl}^{\prime\prime}(P_n)=3$ and
$\chi_g^{\prime\prime}(C_n)=\chi_{gl}^{\prime\prime}(C_n)=4$.
\end{theorem}
\textbf{Proof.} Clearly $\chi_g^{\prime\prime}(P_n)\geq
\chi^{\prime\prime}(P_n)=3$. On the other hand,
$\chi_{gl}^{\prime\prime}(P_n)\leq \mbox{col}(T(P_n))=3.$ It follows
that $\chi_g^{\prime\prime}(P_n)=\chi_{gl}^{\prime\prime}(P_n)=3$.
Denote by $v_1,v_2,\ldots,v_n$ the vertices of $C_n$ in the order
they appear in $C_n$ with $u_i=v_iv_{i+1}$ and $T=T(C_n)$. By
Theorem \ref{Brook's theorem}, we have
$\chi_{gl}^{\prime\prime}(C_n)\leq 4$, for any $n$. Also it is easy
to see \cite{total C_n} that $\chi^{\prime\prime}(C_n)=3$ if $n=3t$
and $\chi^{\prime\prime}(C_n)=4$, otherwise. Consequently, if $n$ is
not a multiple of three, we obtain the desired result. So let $n=3t$
and $A=(Z_3,+)$ be the group with elements 0,1,2 where "+" is the
addition modulo 3. Define $f\in F(T,Z_3)$ with $f(u_{n-1}u_n)=1$,
$f(u_nu_1)=2$ and $f(e)=0$ otherwise, where the edge $u_{n-1}u_n$ in
$T$ is directed from $u_{n-1}$ to $u_{n}$ and the edge $u_nu_1$ is
directed from $u_{n}$ to $u_{1}$. Suppose that
$c:V(T)\longrightarrow Z_3$ is an $(Z_3,f)$-coloring of the vertices
of $T(C_n)$ with $c(v_1)=i$, $c(v_2)=j$ and $c(u_1)=k$, where $0\leq
i\neq j\neq k\leq 2$. Since $f(u_1u_2)=f(v_2u_2)=0$, $c(u_2)$ must
be different from $c(u_1)=k$ and $c(v_2)=j$ and hence $c(u_2)=i$. By
the same reasoning, $c(v_3)=k$ and since $n=3t$ we have
$c(v_{n-1})=j$, $c(v_n)=k$ and finally $c(u_{n-1})=i$. As a
consequence, for any choices of $i,j,k$, since $f(u_{n-1}u_n)=1$ and
$f(u_nu_1)=2$, the vertex $u_n$ can not admit any admissible color
and hence $\chi_{g}^{\prime\prime}(C_n)\geq 4$. Therefore for any
$n$, we obtain $\chi_{g}^{\prime\prime}(C_n)=4$. Now the inequality
$\chi_{g}^{\prime\prime}(C_n)\leq \chi_{gl}^{\prime\prime}(C_n)\leq
4$, implies that
$\chi_g^{\prime\prime}(C_n)=\chi_{gl}^{\prime\prime}(C_n)=4$ which
completes the proof.

$\hfill \blacksquare$

\medskip
\noindent Note that if $G$ is a graph with components $G_1,G_2,\ldots,G_t$, then we have:
\medskip

$\chi_g^{\prime\prime}(G)=\max\{\chi_g^{\prime\prime}(G_1),\ldots,\chi_g^{\prime\prime}(G_t)\}
,\qquad\chi_{gl}^{\prime\prime}(G)=\max\{\chi_{gl}^{\prime\prime}(G_1),\ldots,\chi_{gl}^{\prime\prime}(G_t)\}$.

\medskip \noindent Combining these facts and Theorem \ref{total of C_n}, we obtain the
following corollary.

\medskip
\begin{corollary}\rm{}\label{degree at most two}
Let $G$ be a graph with maximum degree at most two. Then we have
$\chi_g^{\prime\prime}(G)=\chi_{gl}^{\prime\prime}(G)$.
\end{corollary}

\medskip
The {\it wheel graph}, $W_n$, is the graph obtained from $C_n$ by
adjoining a vertex to all vertices of $C_n$. For $n\geq6$, it is
easy to see that $\mbox{col}(T(W_n))=n+1$ and hence
$\chi_g^{\prime\prime}(G)=\chi_{gl}^{\prime\prime}(G)=n+1$. So these
graphs satisfy Conjectures \ref{total conjecture 1} and \ref{total
conjecture 2}. It seems that for every graph $G$ with unique vertex of maximum degree, $\chi_g^{\prime\prime}(G)=\chi_{gl}^{\prime\prime}(G)=\Delta+1$.


\begin{theorem}\rm{}(\cite{Naserasar})\label{naserasr}
For every planar graph $G$ with minimum degree at least 3 there is
an edge $e=uv$ with $\deg(u)+\deg(v)\leq13$.
\end{theorem}


\begin{theorem}\rm{}\label{planar}
Let $k\geq 11$ and $G$ be a planar graph with maximum degree at most
$k$. Then $\mbox{col}(T(G))\leq k+2$.
\end{theorem}
\textbf{Proof.} Let $G$ be a minimal counterexample for Theorem
\ref{planar}. Then there exists a $k\geq 11$ so that
$\mbox{col}(T(G))> k+2$. If $G$ contains a vertex $u$ with
$\deg(u)\leq 2$, then by the minimality of $G$ we have
$\mbox{col}(T(G-u))\leq k+2$. Since the degree of each edge incident
to $u$ in $T(G)-u$ is at most $k+1$ we have $\mbox{col}(T(G))\leq
k+2$, which contradicts the choice of $G$ as a counterexample. So we
may assume that $\delta\geq 3$. By Theorem \ref{naserasr}, there exists
an edge $e=uv$ with $\deg_{T(G)}(e)\leq 13$. We may assume that
$\deg(u)\leq 6$. By minimality of $G$ we have $\mbox{col}(T(G-e))\leq k+2$.
Since $\max\{\deg_{T(G)-u}(e),\deg_{T(G)}(u)\}\leq 12$ and $k\geq 11$
we obtain that $\mbox{col}(T(G))\leq k+2$, which is a contradiction.
 $\hfill \blacksquare$\\

Using Theorem \ref{planar}, we obtain the following corollary, which
states that Conjecture \ref{total conjecture 1} holds for planar
graphs with maximum degree at least 11.

\medskip
\begin{corollary}\rm{}
Let $G$ be a planar graph with maximum degree $\Delta$. Then
$$\chi_{gl}^{\prime\prime}(G)\leq
\mbox{col}(T(G))\leq \max\{13,\Delta+2\}.$$
\end{corollary}
\section{\large{Outerplanar and near outerplanar graphs}}

\vspace{.4 cm}
\hspace{.5 cm} In this section, we give some upper bounds for the
total group choice number of outerplanar and near outerplanar
graphs, which establish Conjectures \ref{total conjecture
1} and \ref{total conjecture 2} for these graphs. We need the
following lemma of Borodin and Woodall \cite{outerplane graphs}.

\begin{lemma}\rm\label{outerplane}
Let $G$ be an outerplanar graph. Then at least one of the following holds.
\\\\(a) $\delta(G)=1$.\vspace*{0.1cm}
\\ (b) There exists an edge $uv$ such that $\deg(u)=\deg(v)=2$. \vspace*{0.1cm}
\\ (c) There exists a 3-face $uxy$ such that $\deg(u)=2$ and $\deg(x)=3$.\vspace*{0.1cm}
\\ (d) There exist two 3-faces $xu_1v_1$ and $xu_2v_2$ such that $\deg(u_1)=\deg(u_2)=2$ and $\deg(x)=4$ and
these five vertices are all distinct.
\end{lemma}

The following theorem implies that Conjectures \ref{total conjecture
1} and \ref{total conjecture 2} hold for outerplanar graphs with
maximum degree at least 5.

\begin{theorem}\rm{}\label{total outerplane}
Let $k\geq 5 $ and $G$ be an outerplanar graph with maximum degree
$\Delta\leq k$. Then $\mbox{col}(T(G))\leq k+1$.
\end{theorem}
\textbf{Proof.} Let $G$ be a minimal counterexample for the theorem.
So $\mbox{col}(T(G))> k+1$ for some $k\geq 5$. If $G$ contains a
vertex $v$ of degree one with neighborhood $u$, then by minimality
of $G$ we have $\mbox{col}(T(G-v))\leq k+1$. Since
$\deg_{T(G)-v}(uv)\leq k$ and $\deg_{T(G)}(v)\leq 2$ we have
$\mbox{col}(T(G))\leq k+1$, a contradiction. So by Lemma
\ref{outerplane}, $G$ contains an edges $uv$ such that $\deg(u)=2$
and $\deg(u)+\deg(v)\leq 6$. By the minimality of $G$,
$\mbox{col}(T(G-uv))\leq k+1$. Again since $k\geq 5$,
$\deg_{T(G)-u}(uv)\leq 5$ and $\deg_{T(G)}(u)=4$, we obtain that
$\mbox{col}(T(G))\leq k+1$, a contradiction. $\hfill \blacksquare$\\

Combining Theorems \ref{total outerplane} and \ref{bound in term of
col(G)}, we have the following corollary.

\medskip
\begin{corollary}\rm{}
Let $G$ be an outerplanar graph with maximum degree $\Delta\neq
4$,
then $\chi_{gl}^{\prime\prime}(G)\leq \max\{5,\Delta+1\}$. In
particular if $\Delta\geq 5$ then
$\chi_{gl}^{\prime\prime}(G)=\chi_{g}^{\prime\prime}(G)=\Delta+1$.
\end{corollary}

\medskip
By a {\it near outerplanar} graph we mean one that is either
$K_4$-minor free or $K_{2,3}$-minor free. Near outeplanar graphs
are an extension of outerplanar graphs. Theorem \ref{near
outreplane}, will show that Conjecture \ref{total conjecture 1}
holds for the class of $K_{2,3}$-minor free graphs. In fact, in
Theorem \ref{near outreplane} we will replace the class of
$K_{2,3}$-minor free graphs by the slightly larger class of
$\bar{K_2}+(K_1\cup K_2)$-minor free graphs, where
$\bar{K_2}+(K_1\cup K_2)$ is the graph obtained from $K_{2,3}$ by
adding an edge joining two vertices of degree 2. Before we proceed, we
need the following lemma which characterizes $\bar{K_2}+(K_1\cup
K_2)$-minor free graphs.

\medskip
\begin{lemma}\rm{}(\cite{near outreplane})\label{near}
Let $G$ be a $\bar{K_2}+(K_1\cup K_2)$-minor free graph. Then each block of $G$ is either $K_4$-minor free or isomorphic to $K_4$.
\end{lemma}

\medskip
\begin{theorem}\rm{}\label{near outreplane}
Let $k\geq 4$ and $G\neq K_4$ be a $\bar{K_2}+(K_1\cup K_2)$-minor
free graph with maximum degree $\Delta\leq k$. Then
$\mbox{col}(T(G))\leq k+2$.
\end{theorem}
\textbf{Proof.} Let $G$ be a minimal counterexample for Theorem
\ref{near outreplane} and also let $k\geq 4$ be such that $\mbox{col}(T(G))>
k+2$. We may assume that $G$ is connected. First, let $G$ be
2-connected. By Lemma \ref{near} we may assume that $G$ is a
$K_4$-minor free graph. So by Corollary \ref{degree at most two},
$G$ is non-regular and hence $\mbox{col}(T(G))\leq \Delta
+\mbox{col}(G)-1\leq k+2$, a contradiction. Hence $G$ is not a
2-connected graph. Suppose $B$ is an end-block of $G$ with
cut-vertex $v$. First let $B\cong K_4$ with $V(B)=\{v,u,w,x\}$. By
minimality of $G$ we have $\mbox{col}(T(G-(B-v)))\leq k+2$. Since
$\deg_{H}(e)\leq k$ for $e\in N_{T(G)}(v)\cap E(B)$ and
$H=T(G)-\{u,w,x,uw,ux,wx\}$, we have $\mbox{col}(T(G))\leq k+2$, a
contradiction. So $B$ is $K_4$-minor free, and this means that $B$
contains at least two vertices of degree at most 2. Let $v$ be a
vertex of degree 2 in $B$ such that $\deg_G(v)=2$. By minimality of
$G$ we have $\mbox{col}(T(G-v))\leq k+2$. Since $\deg_{T(G)-v}(e)\leq k+1$
for $e\in N_{T(G)}(v)$, we have $\mbox{col}(T(G))\leq k+2$, a
contradiction.

$\hfill \blacksquare$

\medskip
By Theorem \ref{D-choosable} we have $\chi_{gl}^{\prime\prime}(K_4)\leq 6$
and so by Theorem \ref{near outreplane}, Conjecture
\ref{total conjecture 1} holds for $K_{2,3}$-minor free graphs with
maximum degree at least 4.

\bigskip
\begin{corollary}\rm{}\label{K_{2,3}}
Let $G$ be a $K_{2,3}$-minor free graph with maximum degree at least
4. Then $G$  is $(\Delta+ 2)$-total group choosable.
\end{corollary}

\bigskip
Using Corollaries \ref{col(G)<3} and  \ref{K_{2,3}} we
conclude that Conjecture \ref{total conjecture 1} holds for near-outerplanar graphs.


\section{\large{Planar graphs without small cycles}}

\bigskip
\hspace{.5 cm}In this section, we prove that Conjecture \ref{total
conjecture 1} holds for some planar graphs without certain small
cycles. Our proofs are based on some structure lemmas and discharging
method.

\medskip
\begin{lemma}\rm{}\label{priciple lemma}
Let $G$ be a graph with $\delta\leq 2$ and $\Delta\geq 3$. If for
any $e\in E(G)$, $\mbox{col}(T(G-e))\leq \Delta+2$, then
$\mbox{col}(T(G))\leq \Delta+2$.
\end{lemma}
\textbf{Proof.} Let $v$ be a vertex of degree at most two and $e$ be
an edge incident with $v$. By the hypothesis,
$\mbox{col}(T(G-e))\leq \Delta+2$. Since $\deg_{T(G)-v}(e)\leq
\Delta+1$ and $\deg_{T(G)}(v)\leq 4$, we have $\mbox{col}(T(G))\leq
\Delta+2$. $\hfill \blacksquare$\\

A {\it plane graph} is a particular drawing of a planar graph on the
plane. A 2{\it-alternating cycle} in a graph $G$ is a cycle of even
length in which alternate vertices have degree 2 in $G$.

\medskip
\begin{theorem}\rm{}(\cite{5-cycles or 6-cycles})\label{structure of 5-cycle or 6 cycle}
Let $G$ be a connected planar graph with $\delta\geq 2$. If $G$
contains no 5-cycles nor 6-cycles, then $G$ contains a 2-alternating
cycle or an edge $uv$ such that $\deg(u)+\deg(v)\leq 9$.
\end{theorem}

\medskip
\noindent
The following theorem establishes Conjecture \ref{total
conjecture 1} for a planar graph with maximum degree at least 7 that
contains no 5-cycles or 6-cycles.

\medskip
\begin{theorem}\rm{}\label{5-cycle or 6 cycle}
Let $k\geq 7$ and $G$ be a planar graph with maximum degree
$\Delta\leq k$. If $G$ contains no 5-cycles or 6-cycles, then
$\mbox{col}(T(G))\leq k+2$.
\end{theorem}
\textbf{Proof.} Let $G=(V,E)$ be a minimal counterexample for
Theorem \ref{5-cycle or 6 cycle}. So  for some $k\geq 7$,
$\mbox{col}(T(G))> k+2$. If $G$ contains a vertex $v$ of degree one,
then by the minimality of $G$, $\mbox{col}(T(G-v))\leq k+2$. Since
$\deg_{T(G)-v}(e)\leq k$ and $\deg_{T(G)}(v)\leq 2$ where $e$ is the
edge incident with $v$, we have $\mbox{col}(T(G))\leq k+2$, a
contradiction. So we may assume that $\delta\geq 2$. First suppose
that $G$ has an edge $e=uv$ with $\deg(u)+\deg(v)\leq 9$. Without
loss of generality, we assume that $\deg(u)\leq 4$. Then
$\mbox{col}(T(G-e))\leq k+2$. Since $\deg_{T(G)-u}(e)\leq 8$ and
$\deg_{T(G)}(u)\leq 8$, we have $\mbox{col}(T(G))\leq k+2$, a
contradiction. Hence for every edge $e=uv$ of $G$,
$\deg(u)+\deg(v)\geq 10$. By Theorem \ref{structure of 5-cycle or 6
cycle}, $G$ must contain a 2-alternating cycle $C$. Let $U$ be the
set of the vertices of $C$ that have degree 2 in $G$, and let
$H=G-U$. By the minimality of $G$ and as $|V(H)|<|V(G)|$, we have
$\mbox{col}(T(H))\leq k+2$. For each $e\in E(C)$ and $v\in U$ of $C$
in $G$ we have $\deg_{T(G)-U}(e)\leq k+1$ and $\deg_{T(G)}(v)\leq 4$
and so $\mbox{col}(T(G))\leq k+2$, a contradiction. This
contradiction completes the proof of the theorem.

$\hfill \blacksquare$

\newpage
A cycle $C$ of length $k$ in a graph $G$ is called a {\it $k$-net}
if $C$ has at least one chord in $G$. The following, is a structural
lemma for plane graphs without $5$-nets.

\begin{lemma}\rm{}(\cite{no 5-nets})\label{no 5-nets}
Let $G$ be a planar graph with $\delta\geq 3$ and without 5-nets. Then $G$ contains an edge $xy$ such that $\deg(x)+\deg(y)\leq 9$.
\end{lemma}
\noindent

The following theorem establishes Conjecture \ref{total conjecture
1} for every planar graph without 5-nets and maximum degree at
least 7.

\begin{theorem}\rm{}\label{5-nets group coloring}
Let $k\geq 7$ and $G$ be a planar graph with maximum degree
$\Delta\leq k$. If $G$ contains no 5-nets, then
$\mbox{col}(T(G))\leq k+2$.
\end{theorem}
\textbf{Proof.} Let $G=(V,E)$ be a minimal counterexample for Theorem
\ref{5-nets group coloring}. Then there is a $k\geq 7$ so that
$\mbox{col}(T(G))> k+2$. By Lemma \ref{priciple lemma}, we may
assume that $\delta\geq 3$. Using Lemma \ref{no 5-nets}, $G$
contains an edge $e=xy$ such that $\deg(x)+\deg(y)\leq 9$ and
$\deg(x)\leq 4$. By the minimality of $G$ we have $\mbox{col}(T(G-e))\leq
k+2$. Since $\deg_{T(G)-x}(e)\leq 8$ and $\deg_{T(G)}(x)\leq 8$, we
have $\mbox{col}(T(G))\leq k+2$, a contradiction.
$\hfill\blacksquare$\\

We denote the set of faces of a plane graph $G$ by $F(G)$ or simply
by $F$. For a plane graph $G$ and $f \in F(G)$, we write
$f=u_1u_2\ldots u_n$ if $u_1, u_2, \ldots, u_n$ are the vertices on
the boundary walk of $f$ enumerated clockwise. Let $\delta(f)$
denote the minimum degree of vertices incident with $f$. The
{\it{degree of a face}} $f$, denoted by $\deg(f)$, is the number of
edge steps in the boundary walk. A {\it $k$-vertex} (resp. {\it $k^+$-vertex})
is a vertex of degree $k$ (resp. a vertex of degree at least $k$).
The following theorem establishes Conjecture \ref{total conjecture
1} for planar graphs without 4-cycles and maximum degree at
least 6.

\begin{theorem}\rm{}\label{no 4-cycles group coloring}
Let $k\geq 6$ and $G$ be a planar graph with maximum degree
$\Delta\leq k$ such that $G$ has no cycle of length 4. Then
$\chi_{gl}^{\prime\prime}(G)\leq k+2$.
\end{theorem}
\textbf{Proof.} Let $G$ be a minimal counterexample. So for a $k\geq
6$, a group $A$ with $|A|\geq k+2$, a $(k+2)$-list assignment
$L:V(T(G))\longrightarrow 2^A$ and $f\in F(T(G),A)$, $T(G)$ is not
$(A,L,f)$-colorable. Graph $G$ has the following properties:
\noindent
\\(1) $G$ is connected,
\\(2) Any vertex $v$ is incident with at most $\lfloor\frac{\deg(v)}{2}\rfloor$ 3-faces,
\\(3) The minimum degree of $G$ is at least three that is $\delta\geq 3$,
\\(4) $G$ contains no edge $uv$ with $\min\{\deg(u),\deg(v)\}\leq \frac{k}{2}$ and $\deg(u)+\deg(v)\leq k+2$,
\\(5) $G$ does not contain any 3-face $F=uvw$ such that $\deg(u)=\deg(v)=\deg(w)=4$.

\medskip The proofs of (1) and (2) are clear. If $G$ contains a vertex $v$ of
degree at most two and $e$ is an edge incident with $v$, then
$\chi_{gl}^{\prime\prime}(G-e)\leq k+2$, by minimality of $G$, and
so there exists an $(A,L,f)$-coloring $c$ for $T(G-e)$. Erase the
color of vertex $v$ in this coloring. There are at least
$(k+2)-(k+1)$ usable colors in the list of edge $e$, and so it can
be colored. Now since $k\geq 6$, the vertex $v$ can be colored from
its list, and this gives an $(A,L,f)$-coloring for $T(G)$, which
contradicts the minimality of $G$.  Also if $G$ contains an edge
$uv$ with $\min\{\deg(u),\deg(v)\}\leq \frac{k}{2}$ and
$\deg(u)+\deg(v)\leq k+2$, then any $(A,L,f)$-coloring of $T(G-uv)$
can be extended to an $(A,L,f)$-coloring of $T(G)$, which is
contradiction. This shows that (3) and (4) hold. To see $(5)$, on
the contrary, let such a face exists. Let
$G^{\prime}=G-\{uv,vw,uw\}$. By the minimality of $G$,
$T(G^{\prime})$ has an $(A, L, f)$-coloring $c$. Erase the colors of
$u, v, w$ and, for an element $x \in \{u, v, w, uv, vw, uw\}$, let
$L'(x)$ be the available colors in the list of $x$. Since $k \ge 6$
and since $u, v$, and $w$ are degree 4 vertices, each $|L'(x)| \ge k
- 2 \ge 4$. By Theorem \ref{total of C_n}, it is possible to recolor
the elements $u, v, w, uv, vw, uw,$ in this path with colors $c(u)
\in L'(u), c(v) \in L'(v)$, $c(w)\in L'(w)$, $c(uv) \in L'(uv)$,
$c(vw) \in L'(vw)$ and $c(uw) \in L'(uw)$, respectively, so that $c$
is indeed an $(A, L, f)$-coloring of $T(G)$, which is a
contradiction.

Since $G$ is a planar graph, by Euler's Formula, we have:

\bigskip
$\sum_{v\in V}(2\deg(v)-6)+\sum_{f\in F}(\deg(f)-6)=-6(|V
|-|E|+|F|)=-12.$

\bigskip\noindent We define the initial charge function $w(x)$ for
each $x\in V\cup F$. Let $w(v)=2\deg(v)-6$ if $v\in V$ and
$w(f)=\deg(f)-6$ if $f\in F$. It follows that $\sum _{_{x\in V\cup
F}}w(x)< 0$. We construct a new charge $w^*(x)$ on $G$ as follows:
\medskip
\\Each 3-face receives $\frac{3}{2}$ from its incident vertices of
degree at least 5.
\\Each 3-face receives $\frac{3}{4}$ from its incident vertices of
degree 4.
\\Each 5-face receives $\frac{1}{3}$ from its incident vertices of
degree at least 5.
\\Each 5-face receives $\frac{1}{4}$ from its incident vertices of
degree 4.

\medskip
Note that $w^*(f )=w(f )\geq 0$ if $\deg(f)\geq6$. Assume that
$\deg(f)=3$. If $\delta(f)=3$, then $f$ is incident with two
$6^+$-vertices by (4). So $w^*(f)\geq w(f)+2\times\frac{3}{2}=0$.
Otherwise, $f$ is incident with at least one $5^+$-vertex by (5). So
$w^*(f)\geq w(f)+\frac{3}{2}+2\times\frac{3}{4}=0$. Let $\deg(f)=5$.
If $\delta(f)=3$, then $f$ is incident with at most two vertices of
degree 3 by (4), and if $f$ is incident with two vertices of degree
3, then $f$ is incident with three $6^+$-vertices. Thus $w^*(f)\geq
w(f)+\min\{2\times\frac{1}{3}+2\times
\frac{1}{4},3\times\frac{1}{3}\}=0$. Otherwise, $w^*(f)\geq w(f)+
5\times\frac{1}{4}>0$. Let $v$ be a vertex of $G$. Clearly,
$w^*(v)=w(v)=0$ if $\deg(v)=3$. If $\deg(v)=4$, then $v$ is incident
with at most two 3-faces by (2). So $w^*(v)\geq
w(v)-2\times\frac{3}{4}-2\times\frac{1}{4}=0$. If $\deg(v)=5$, then
$v$ is incident with at most two 3-faces by (b). So $w^*(v)\geq
w(v)-2\times\frac{3}{2}-3\times\frac{1}{3}=0$. If $\deg(v)=6$, then
$w^*(v)\geq w(v)-3\times\frac{3}{2}-3\times\frac{1}{3}>0$. If
$\deg(v)\geq7$, so $w^*(v)\geq
w(v)-\lfloor\frac{\deg(v)}{2}\rfloor\times\frac{3}{2}-\lceil
\frac{\deg(v)}{2}\rceil\times\frac{1}{3}>0.$ It follows that $\sum
_{_{x\in V\cup F}}w^*(x)=\sum _{_{x\in V\cup F}}w(x)\geq 0$, a
contradiction. $\hfill \blacksquare$\\

A face $f$ is called {\it simple} if its boundary is a cycle. If
$f=u_1u_2\ldots u_n$ is not simple, then $f$ contains at least one
cut vertex $v$. Let $m_v(f)$ denotes the number of times through $v$
of $f$ in clockwise order. In the sequel, we prove that plane graphs
without 4,5-cycles with maximum degree $\Delta\geq 5$ are totally
$(\Delta+2)$-group choosable.

\medskip
\begin{theorem}\rm{}
Let $k\geq 5$ and $G$ be a planar graph with maximum degree
$\Delta\leq k$ such that $G$ has no cycles of length 4 and 5. Then
$\chi_{gl}^{\prime\prime}(G)\leq k+2$.
\end{theorem}
\textbf{Proof.} Let $G$ be a minimal counterexample. For a $k\geq
5$, a group $A$ with $|A|\geq k+2$, a $(k+2)$-list assignment
$L:V(T(G))\longrightarrow 2^A$ and $f\in F(T(G),A)$, $T(G)$ is not
$(A,L,f)$-colorable. Theorem \ref{no 4-cycles group coloring}
implies the theorem when $k \geq 6$. Hence it is sufficient to prove
the theorem when $k=5$. Graph $G$ has the following properties:
\medskip
\noindent
\\(1) $G$ is connected,
\\(2) Any vertex $v$ is incident with at most $\lfloor\frac{\deg(v)}{2}\rfloor$ 3-faces,
\\(3) The minimum degree of $G$ is at least 3 i.e $\delta\geq 3$,
\\(4) $G$ contains no edge $uv$ with $\min\{\deg(u),\deg(v)\}=3$ and $\deg(u)+\deg(v)\leq 7$.

\medskip The proofs of (1)-(4) are similar to the proof of Theorem \ref{no
4-cycles group coloring}. We define the initial charge function
$w(x)=\deg(x)-4$ for each $x\in V\cup F$. By the Euler's Formula, we
have $\sum _{_{x\in V\cup F}}w(x)=-8$. We construct a new charge
$w^*(x)$ on $G$ as follows:\\
Each $r(\geq 6$)-face $f$ gives $(1-\frac{4}{r})m_v(f)$ to its
incident vertex $v$ if $v$ is cut vertex, and gives $1-\frac{4}{r}$
otherwise.\\
Each 3-vertex $v$ receives $\frac{1}{3}$ from $u$ if $v$ is incident
with 3-face $f$ and $u$ is a neighbor of $v$ but not incident with
$f$.\\
Each 3-face receives $\frac{1}{2}$ from its incident vertex $v$ if
$\deg(v)=5$ and receives $\frac{1}{3}$ if $\deg(v)=4$.\vspace*{.2cm}

 Note that $w^*(f)\geq 0$ for any face $f$. Let $v$ be a
vertex of $G$. Suppose that $\deg(v)=3$. If $v$ is incident with a 3-face
$f$, then $v$ receives at least $\frac{2}{3}$ from its incident
faces and $\frac{1}{3}$ from its incident vertex not lying on $f$.
So $w^*(v)\geq w(v)+\frac{2}{3}+\frac{1}{3}=0$. Otherwise, $v$
receives at least $3\times\frac{1}{3}$ from its incident faces and
hence $w^*(v)\geq w(v)+1=0$. Let $\deg(v)=4$. The vertex $v$ is
incident with at most two 3-faces by (2), so $v$ gives at most
$\frac{2}{3}$ to its incident 3-faces. Also $v$ receives at least
$\frac{2}{3}$ from its incident faces of degree at least 6. Hence
$w^*(v)\geq w(v)+\frac{2}{3}-\frac{2}{3}=0$. Finally let
$\deg(v)=5$. The vertex $v$ is incident with at most two 3-faces by
(2), and if $v$ is incident with exactly two 3-faces and a 3-vertex
is pending on the remaining neighbors of $v$, then $v$ gives at most
$2\times\frac{1}{2}+\frac{1}{3}$ and receives at least $3\times
\frac{1}{3}$ from its faces of degree at least 6. So $w^*(v)\geq
w(v)+3\times \frac{1}{3}-(2\times\frac{1}{2}+\frac{1}{3})>0$. If $v$
is incident with a 3-face and three 3-vertices is pending on the
remaining neighbors of $v$, then $v$ gives at most
$3\times\frac{1}{3}+\frac{1}{2}$ and receives at least $4\times
\frac{1}{3}$ from its faces of degree at least 6. So $w^*(v)\geq
w(v)+4\times \frac{1}{3}-(3\times\frac{1}{3}+\frac{1}{2})>0$.
Otherwise $v$ gives at most $5\times\frac{1}{3}$ and receives at
least $5\times \frac{1}{3}$ from its faces of degree at least 6. So
$w^*(v)\geq w(v)+\frac{5}{3}-\frac{5}{3}>0$. It follows that $\sum
_{_{x\in V\cup F}}w^*(x)=\sum _{_{x\in V\cup F}}w(x)>0$, a
contradiction. This contradiction completes the proof of the theorem.$\hfill
\blacksquare$

\begin{lemma}\rm{}(\cite{no five cycles})\label{no five cycles}
Let $G$ be a planar graph with $\delta\geq 3$ and no five cycles.
Then there exists an edge $xy$ such that $\deg(x)=3$ and
$\deg(y)\leq 5$.
\end{lemma}

The following theorem proves Conjecture \ref{total conjecture
1} for planar graphs without 5-cycles and maximum degree at
least 6.

\begin{theorem}\rm{}\label{no 5-cycles group coloring}
Let $k\geq 6$ and $G$ be a planar graph with maximum degree
$\Delta\leq k$. If $G$ contains no 5-cycles, then
$\chi_{gl}^{\prime\prime}(G)\leq k+2$.
\end{theorem}
\textbf{Proof.} Let $G=(V,E)$ be a minimal counterexample for
theorem. So for a $k\geq 6$, a group $A$ with $|A|\geq k+2$, a
$(k+2)$-list assignment $L:V(T(G))\longrightarrow 2^A$ and $f\in
F(T(G),A)$, $T(G)$ is not $(A,L,f)$ colorable. By Lemma
\ref{priciple lemma}, we may assume that $\delta\geq 3$. So by Lemma
\ref{no five cycles}, $G$ contains an edge $xy$ such that
$\deg(x)=3$ and $\deg(y)\leq 5$. the graph $G-e$ has an
$(A,L,f)$-coloring by the minimality of $G$. Now erase the color of
$x$ in this coloring and color the edge $xy$ from its list, which is
possible since its list has at least $(k+2)-7\geq 1$ usable colors.
Since $\deg_{T(G)}(x)=6$ and $k\geq 6$, the vertex $x$ can be
colored easily. So $G$ has an $(A,L,f)$-coloring, which is a
contradiction.

$\hfill \blacksquare$

\vspace{1.5cm}

 \end{document}